# $L_{1/2}$ Regularization: Convergence of Iterative Half Thresholding Algorithm

Jinshan Zeng, Shaobo Lin, Yao Wang, and Zongben Xu

*Abstract*—In recent studies on sparse modeling, the nonconvex regularization approaches (particularly, $L_q$ regularization with $q \in (0,1)$) have been demonstrated to possess capability of gaining much benefit in sparsity-inducing and efficiency. As compared with the convex regularization approaches (say, $L_1$ regularization), however, the convergence issue of the corresponding algorithms are more difficult to tackle. In this paper, we deal with this difficult issue for a specific but typical nonconvex regularization scheme, the $L_{1/2}$ regularization, which has been successfully used to many applications. More specifically, we study the convergence of the iterative *half* thresholding algorithm (the *half* algorithm for short), one of the most efficient and important algorithms for solution to the $L_{1/2}$ regularization. As the main result, we show that under certain conditions, the *half* algorithm converges to a local minimizer of the $L_{1/2}$ regularization, with an eventually linear convergence rate. The established result provides a theoretical guarantee for a wide range of applications of the *half* algorithm. We provide also a set of simulations to support the correctness of theoretical assertions and compare the time efficiency of the *half* algorithm with other known typical algorithms for $L_{1/2}$ regularization like the iteratively reweighted least squares (IRLS) algorithm and the iteratively reweighted $l_1$ minimization (IRL1) algorithm.

*Index Terms*—Convergence, iterative *half* thresholding algorithm, $L_q/L_{1/2}$ regularization, nonconvex regularization.

## I. INTRODUCTION

THE sparsity problems emerging in many areas of scientific research and engineering practice have attracted considerable attention in recent years ([1]–[7]). In a typical setup, an unknown sparse vector $x \in \mathbf{R}^N$ is reconstructed from measurements

$$y = Ax, \qquad (1)$$

Manuscript received July 27, 2013; revised November 01, 2013, January 26, 2014, and February 14, 2014; accepted February 19, 2014. Date of publication February 28, 2014; date of current version April 07, 2014. The associate editor coordinating the review of this manuscript and approving it for publication was Prof. Min Dong. The research was supported by the National 973 Programs (Grant No. 2013CB329404, 2010CB731905), the Key Program of National Natural Science Foundation of China (Grants No. 11131006), and the National Natural Science Foundations of China (Grants No. 61075054, 11001227, 11171272). This work was partially supported by Southwest China Institute of Electronic Techonology.

J. Zeng is with the Institute for Information and System Sciences, School of Mathematics and Statistics, Xi'an Jiaotong University, Xi'an 710049, and also with the Beijing Center for Mathematics and Information Interdisciplinary Sciences (BCMIIS), Beijing, 100048, China (e-mail: jsh.zeng@gmail.com).

S. Lin, Y. Wang, and Z. Xu are with the Institute for Information and System Sciences, School of Mathematics and Statistics, Xi'an Jiaotong University, Xi'an 710049, China (e-mail: sblin1983@gmail.com; yao.s.wang@gmail.com; zbxu@mail.xjtu.edu.cn).

Digital Object Identifier 10.1109/TSP.2014.2309076

or more generally, from

$$y = Ax + \epsilon, \qquad (2)$$

where $y \in \mathbf{R}^m$, $A \in \mathbf{R}^{m \times N}$ (commonly, $m < N$) is a measurement matrix and $\epsilon$ represents the noise. The problem can be modeled as the so called $L_0$ regularization problem

$$\min_{x \in \mathbf{R}^N} \left\{ \|Ax - y\|_2^2 + \lambda \|x\|_0 \right\}, \qquad (3)$$

where $\|x\|_0$, formally called $L_0$ quasi-norm, denotes the number of nonzero components of $x$, and $\lambda > 0$ is a regularization parameter. In general, $L_0$ regularization is challenging to solve, due to it is NP-hard [8]. In order to overcome such difficulty, $L_1$ regularization was proposed as an alternative ([1], [2], [9], [10])

$$\min_{x \in \mathbf{R}^N} \left\{ \|Ax - y\|_2^2 + \lambda \|x\|_1 \right\}, \qquad (4)$$

where $\|x\|_1 = \sum_{i=1}^{N} |x_i|$ represents the $L_1$ norm.

It is well known that $L_1$ regularization is a convex optimization problem, and therefore, can be solved very efficiently. It can also result in sparse solution of the considered problem, with the guarantee that, under certain conditions, the resultant solution coincides with the solution of $L_0$ regularization ($L_1/L_0$ equivalence) [11]. Because of these, $L_1$ regularization gets its popularity and has been widely used in various applications such as signal/image processing, biomedical informatics and computer vision [12]–[14]. Nevertheless, $L_1$ regularization has been shown to be suboptimal in many cases (particularly, it cannot enforce further sparsity when applied to compressed sensing) ([15]–[20]), since the $L_1$ norm is a loose approximation of the $L_0$ norm and often leads to an over-penalized problem. Consequently, some further improvements are required. Among such efforts, a very natural improvement is the suggestion of the use of $L_q$ ($0 < q < 1$) regularization ([17]–[20])

$$\min_{x \in \mathbf{R}^N} \left\{ \|Ax - y\|_2^2 + \lambda \|x\|_q^q \right\}, \qquad (5)$$

where $\|x\|_q$ represents the $L_q$ quasi-norm, defined by $\|x\|_q = \left( \sum_{i=1}^{N} |x_i|^q \right)^{1/q}$.

The $L_q$ regularization is a nonconvex and nonsmooth optimization problem, which is in general difficult to be solved fast and efficiently. Moreover, even when solvable, which $q$ should be selected to yield the best result is also a problem. Recent studies have provided some insights into these problems ([19]–[23]). On one hand, the high efficiency of $L_{1/2}$ regularization and $L_{2/3}$ regularization when applied to image deconvolution was demonstrated in [22] and [23], and, on the other hand, a representative role of $L_{1/2}$ regularization among all the $L_q$ regularization with $q$ in (0,1) was revealed in [21] by a phase diagram study. Furthermore, through developing a thresholding representation theory, an iterative *half* thresholding algorithm (called *half* algorithm in brief) was proposed in [20] for fast solution to $L_{1/2}$ regularization. Inspired by the well developed theoretical properties and the fast algorithm, $L_{1/2}$ regularization has been successfully used to many applications including hyperspectral unmixing [24], synthetic aperture radar (SAR) imaging [25], [26], machine learning [27], [28], gene selection [29] and practical engineering [30].

Although widely applied, the theoretical behavior on convergence of $L_{1/2}$ regularization is still unclear. This is not surprising because strikingly different from the convex case, the convergence analysis of the nonconvex case is more difficult to tackle. The convergence issue, however, is very crucial as it characterizes the circumstance in which an algorithm can or cannot be successfully applied, and thus underlies the feasibility of an approach. By applying a dynamic system methodology, Xu *et al.* [20] has ever verified the convergence of the *half* algorithm to a stationary point under some restricted conditions, but it is not answered the question what the convergent stationary point is. In particular, the following two questions are open:
1) Does the *half* algorithm converge to a local or global minimizer of $L_{1/2}$ regularization?
2) How fast does the *half* algorithm converge?

Our purpose of the present paper is to provide answers to these questions. More specifically, due to the nonconvex feature, it is hardly possible to identify the global minimizer. Instead, we focus on characterizing the conditions under which the algorithm converges to a local minimizer, and how fast the algorithm converges eventually.

As the main contributions of the present research, we first extend the convergence result established in [20] to a more general case, that is, we show that the *half* algorithm converges to a stationary point under the same condition as that of convergence of the iterative *soft* thresholding algorithm (the *soft* algorithm in brief) in $L_1$ regularization case [31]. Thus, the well known convergence result of the convex $L_1$ regularization case has been generalized to the nonconvex $L_{1/2}$ regularization case. We then formulate two sufficient conditions under which the convergent stationary point is a local minimizer of $L_{1/2}$ regularization. The first condition states that the regularization parameter $\lambda$ should be relatively small, which means that the square of $l_2$ term, $\|Ax - y\|_2^2$ should be penalized not too heavily. The second condition concerns that the measurement matrix $A$ satisfies a certain concentration assumption. Particularly, when $A$ satisfies a certain restricted isometry property (RIP), the second condition can be satisfied naturally. Under the same conditions, we further verify that the eventual convergence speed of the *half* algorithm is linear, that is, when the iteration number is sufficiently large, the iterations in the remaining steps converge with at least a linear rate. Thus, under the condition, the *half* algorithm can converge to a local minimizer of $L_{1/2}$ regularization with an eventually linear speed. Finally, we provide simulations to support the correctness of theoretical assertions and compare the convergence speed of the *half* algorithm with other known typical algorithms for $L_{1/2}$ regularization like the iteratively reweighted least squares (IRLS) algorithm and the iteratively reweighted $l_1$ minimization (IRL1) algorithm.

The rest of the paper is organized as follows. In Section II, we introduce the related $L_{1/2}$ regularization theory and the *half* algorithm. In Section III, we state the main convergence results with some remarks. In Section IV, we present an eventually linear convergence rate estimation of the *half* algorithm. In Section V, we conduct the simulations to substantiate the theoretical results in Sections III and IV. In Section VI, we compare the obtained theoretical results with some other known results. We then present the proofs of all the theoretical results in Section VII and conclude the research in Section VIII.

## II. ITERATIVE HALF THRESHOLDING ALGORITHM

In this section, we review the related thresholding representation theory of $L_{1/2}$ regularization and the corresponding *half* algorithm, which then serves as the basis of further analysis in the next sections.

### A. Thresholding Representation Theory

The $L_{1/2}$ regularization is concerned with the nonconvex and nonsmooth optimization model

$$\min_x \left\{ \|Ax - y\|_2^2 + \lambda \|x\|_{1/2}^{1/2} \right\}, \tag{6}$$

where $\lambda > 0$ is a regularization parameter. It was shown in [20] that $L_{1/2}$ regularization permits a thresholding representation for its solution, detailed as the following.

*Lemma 1. ([20]):* With any real parameter $\mu \in (0, \|A\|_2^{-2})$, any minimizer $x^*$ of $L_{1/2}$ regularization satisfies

$$x^* = H_{\lambda\mu, \frac{1}{2}}(B_\mu(x^*)), \tag{7}$$

where $H_{\lambda\mu, \frac{1}{2}}$ is the *half* thresholding operator defined by (8)–(10), $B_\mu(x^*) = x^* - \mu A^T(Ax^* - y)$, and $A^T$ is the transpose of $A$.

In (7), the *half* thresholding operator $H_{\lambda\mu, \frac{1}{2}}$ is defined by

$$H_{\lambda\mu, \frac{1}{2}}(z) = \left(h_{\lambda\mu, \frac{1}{2}}(z_1), \cdots, h_{\lambda\mu, \frac{1}{2}}(z_N)\right)^T, \forall z \in \mathbf{R}^N \tag{8}$$

where

$$h_{\lambda\mu, \frac{1}{2}}(z_i) = \begin{cases} f_{\lambda\mu, \frac{1}{2}}(z_i), & |z_i| > \frac{\sqrt[3]{54}}{4}(\lambda\mu)^{\frac{2}{3}} \\ 0, & otherwise \end{cases} \tag{9}$$

with

$$f_{\lambda\mu, \frac{1}{2}}(z_i) = \frac{2}{3} z_i \left(1 + \cos\left(\frac{2\pi}{3} - \frac{2}{3}\phi_{\lambda\mu}(z_i)\right)\right), \tag{10}$$

and $\phi_{\lambda\mu}(z_i) = \arccos\left(\frac{\lambda\mu}{8}\left(\frac{|z_i|}{3}\right)^{-3/2}\right)$ for $i = 1, 2, \cdots, N$. For the convenience of expression, we also define

$$F_{\lambda,\mu}(\cdot) = H_{\lambda\mu, \frac{1}{2}}(B_\mu(\cdot)).$$

The following two propositions summarize some useful properties of the *half* thresholding operator.

*Proposition 1. ([20]):* For any $u \in \mathbf{R}^1$, $f_{\lambda\mu,1/2}(u)$ has the following properties:
(1) $f_{\lambda\mu,1/2}(u)$ is an odd function.
(2) For any $u > \frac{\sqrt[3]{54}}{4}(\lambda\mu)^{2/3}$, $f_{\lambda\mu,1/2}(u)$ is strictly increasing and $f_{\lambda\mu,1/2}(u) > \left(\frac{\lambda\mu}{2}\right)^{2/3}$.
(3) For any $|u| > \frac{\sqrt[3]{54}}{4}(\lambda\mu)^{2/3}$, $f_{\lambda\mu,1/2}(u)$ solves the equation

$$f_{\lambda\mu,\frac{1}{2}}(u) + \lambda\mu \frac{sign\left(f_{\lambda\mu,\frac{1}{2}}(u)\right)}{4\sqrt{\left|f_{\lambda\mu,\frac{1}{2}}(u)\right|}} = u,$$

where $sign(\cdot)$ is the sign function.

*Proposition 2:* For any fixed point $x^*$ of $F_{\lambda,\mu}$ (that is, $x^* = F_{\lambda,\mu}(x^*)$), let $I$ be the support set of $x^*$, i.e., $I = supp(x^*)$, then there holds

$$\left|A_i^T(Ax^* - y)\right| \leq \frac{\sqrt[3]{54}}{4}\lambda^{\frac{2}{3}}\mu^{-1/3}, \ for \ i \in I^c$$

$$A_i^T(Ax^* - y) = -\frac{\lambda sign(x_i^*)}{4\sqrt{|x_i^*|}}, \ for \ i \in I,$$

where $A_i \in \mathbf{R}^m$ denotes the $i$-th column of $A$ and $A_i^T$ denotes the transpose of $A_i$.

Proposition 2 can be easily derived from the form of the *half* thresholding representation (7)–(10) and Proposition 1(3).

*Remark 1:* For comparison, we list the *hard* thresholding operator for $L_0$ regularization and the *soft* thresholding operator for $L_1$ regularization also as follows:
(a) the *hard* thresholding operator ([32])

$$h_{\lambda\mu,0}(u) = \begin{cases} u, & |u| > \sqrt{\lambda\mu} \\ 0, & otherwise, \end{cases}$$

(b) the *soft* thresholding operator ([33])

$$h_{\lambda\mu,1}(u) = \begin{cases} sign(u)(|u| - \lambda\mu/2), & |u| > \lambda\mu/2 \\ 0, & otherwise. \end{cases}$$

### B. Half Algorithm

From the thresholding representation form (7), an iterative thresholding algorithm, called the *half* algorithm, for solution to $L_{1/2}$ regularization was naturally developed in [20], which reads as

$$x^{(n+1)} = H_{\lambda\mu,\frac{1}{2}}\left(x^{(n)} - \mu A^T\left(Ax^{(n)} - y\right)\right), \quad (11)$$

where $\mu > 0$ is a step size parameter. Note that if, instead of using the *half* thresholding operator by the *soft* and *hard* thresholding operators, (11) then corresponds to the iterative *soft* and *hard* thresholding algorithms for $L_1$ regularization [33] and $L_0$ regularization [32].

The convergence of the *half* algorithm has been partially analyzed in [20], which implies that the *half* algorithm can converge to a stationary point of (7) when $\mu$ is sufficiently small. In the following, we refine this result by relaxing the restriction on $\mu$ to a very loose condition, i.e., $0 < \mu < \|A\|_2^{-2}$, which is the condition widely imposed for convergence of other convex regularization, say, $L_1$ regularization [31].

*Lemma 2:* Let $\{x^{(n)}\}$ be the sequence generated by the *half* algorithm and suppose that the step size $\mu$ satisfies $0 < \mu < \|A\|_2^{-2}$. Then $x^{(n)}$ converges to a stationary point of (7).

The proof of Lemma 2 will be presented in Section VII.

In the next, we present an important property of the *half* algorithm, that is, the support sets of the iterations of (11) will remain stable within finite steps.

*Lemma 3:* Assume that $\{x^{(n)}\}$ is the sequence generated by the *half* algorithm which converges to $x^*$. If $I$ and $I_n$ are respectively the support sets of $x^*$ and $x^{(n)}$, then there exists an $n_0$ such that whenever $n > n_0$, there holds
(1) $\|x^{(n)} - x^*\|_2 < \left(\frac{\lambda\mu}{2}\right)^{2/3}$;
(2) $I_n = I$; and
(3) $sign(x^{(n)}) = sign(x^*)$.

We present the proof of Lemma 3 in Section VII.

## III. CONVERGENCE TO A LOCAL MINIMIZER

In the last section we have justified the convergence of the *half* algorithm to a stationary point of the iteration (7). However, in practice, we may be more interested in whether the algorithm converges to a global minimizer or more possibly, a local minimizer due to the nonconvexity of the regularization problem. In this section, we report the convergence of the algorithm to a local minimizer of $L_{1/2}$ regularization under certain conditions, followed then with some remarks.

### A. Results Characterized With Regularization Parameter

Let $T_\lambda(x) = \|Ax - y\|_2^2 + \lambda\|x\|_{1/2}^{1/2}$. For any vector $x \in \mathbf{R}^N$, we denote by $supp(x)$ the support set of $x$, i.e., $supp(x) = \{i : |x_i| \neq 0, i = 1, \cdots, N\}$. For any index set $S$, we denote by $A_S$ the submatrix of $A$ with the columns restricted to $S$ and $\sigma_{\min}(A_S^T A_S)$ the minimal singular value of $A_S^T A_S$. With those notations, we can state our main result in the following.

*Theorem 1:* Assume that $0 < \mu < \|A\|_2^{-2}$ and the sequence $\{x^{(n)}\}$ generated by the *half* algorithm converges to $x^*$. Let $I = supp(x^*)$ and $e = \min_{i \in I} |x_i^*|$. If $\sigma_{\min}(A_I^T A_I) > 0$ and

$$0 < \lambda < 8e^{\frac{3}{2}}\sigma_{\min}(A_I^T A_I), \quad (12)$$

then there exists a constant $\beta_0 \in (0, 1)$ such that for any $x$ satisfying $\|x - x^*\|_2 < \beta_0 e$, there holds the estimation

$$T_\lambda(x) - T_\lambda(x^*) \geq 0,$$

that is, $x^*$ is a local minimizer of $L_{1/2}$ regularization.

We provide the proof of Theorem 1 in Section VII. Theorem 1 implies that when the regularization parameter $\lambda$ is taken to be relatively small, the convergent stationary point $x^*$ of the *half* algorithm is indeed a local minimizer of $L_{1/2}$ regularization.

Intuitively, a small $\lambda$ means that the $l_2$ term $\|y - Ax\|_2^2$ of $L_{1/2}$ regularization is penalized not too heavily. More specifically, we observe from Lemma 3 that when $\|x - x^*\|_2 < \beta_0 e$ with $\beta_0 \in (0, 1)$, there holds the identity

$$sign(x_i) = sign(x_i^*)$$

for any $i \in I$ and

$$\|x_{I^c}\|_2 < \beta_0 e,$$

where $I^c$ represents the complementary set of $I$. For any $x$ satisfying $\|x - x^*\|_2 < \beta_0 e$, we reformulate $T_\lambda(x)$ as

$$T_\lambda(x) = \left( \|A_{I^c} x_{I^c}\|_2^2 + \lambda \|x_{I^c}\|_{\frac{1}{2}}^{\frac{1}{2}} - 2(y - A_I x_I)^T A_{I^c} x_{I^c} \right)$$
$$+ \left( \|y - A_I x_I\|_2^2 + \lambda \|x_I\|_{\frac{1}{2}}^{\frac{1}{2}} \right)$$
$$:= T_{I^c} + T_I. \quad (13)$$

where $T_{I^c}$ and $T_I$ represent respectively the first and second parts of the above expression of $T_\lambda(x)$. Then we have

$$T_{I^c} = 2(A_I x_I - A_I x_I^*)^T A_{I^c} x_{I^c} + 2(A_I x_I^* - y)^T A_{I^c} x_{I^c}$$
$$+ \|A_{I^c} x_{I^c}\|_2^2 + \lambda \|x_{I^c}\|_{\frac{1}{2}}^{\frac{1}{2}}$$
$$\leq 2\beta_0 e \|A_I\|_2 \|A_{I^c} x_{I^c}\|_2 + 2\|Ax^* - y\|_2 \|A_{I^c} x_{I^c}\|_2$$
$$+ \|A_{I^c} x_{I^c}\|_2^2 + \lambda \|x_{I^c}\|_{\frac{1}{2}}^{\frac{1}{2}},$$

where the inequality was deduced from $\|x_I - x_I^*\|_2 \leq \|x - x^*\|_2 < \beta_0 e$ and $A_I x_I^* = Ax^*$. Furthermore, since $\|x_{I^c}\|_2 < \beta_0 e$ that can be sufficiently small when $\beta_0$ is very small, $T_{I^c}$ then can be very small. Thus, in this case, the second part of (13), $T_I$, becomes dominant in (13). In consequence, to guarantee $x^*$ being a local minimizer of $T_\lambda(x)$, we only need to justify the strict convexity of $T_I$ at $x^*$.

Since the regularization term $\|x_I\|_{1/2}^{1/2}$ is differential at the neighborhood of $x^*$, $T_I$ is strictly convex at $x^*$ if the following second order condition holds, i.e.,

$$2A_I^T A_I - \lambda \Lambda_I^* \succ 0, \quad (14)$$

where $\Lambda_I^*$ is a diagonal matrix with $\Lambda_I^*(i,i) = \frac{sign(x_i^*)}{4|x_i^*|^{3/2}}$ for $i \in I$ and, for any matrix $E$, the notation $E \succ 0$ represents the positive definetness of $E$. Apparently, if $\lambda$ satisfies the condition in Theorem 1, the second order condition (14) is satisfied. Therefore, $x^*$ is a local minimizer of $T_\lambda(x)$.

Moreover, we can observe from (14) that a weaker condition than (12) can be adopted, say,

$$0 < \lambda < \min_{i \in I} \left\{ 8|x_i^*|^{\frac{3}{2}} \sigma_i \left( A_I^T A_I \right) \right\},$$

where $\sigma_i(A_I^T A_I)$ is the $i$-th singular value of $A_I^T A_I$.

*Remark 2:* From the procedure of proof of Theorem 1, we can observe that when $\lambda$ is taken smaller, $\beta_0$ can be larger, which then implies that the attraction basin of $x^*$ gets larger, thus $x^*$ becomes a more stable local minimizer. This fact is important in practice, since it is hardly possible in general to achieve a global minimizer for a nonconvex method, and, instead, a more stable local minimizer is usually expected in application.

### B. Results Characterized With Measurement Matrix

By analyzing the condition (12) in Theorem 1, we can see that the constraints on $\lambda$ is mainly controlled by two quantities, the absolute value of the minimal nonzero component of $x^*$ and the minimal singular value of $A_I^T A_I$. Moreover, by Proposition 1(2), we have

$$e > \left( \frac{\lambda \mu}{2} \right)^{\frac{2}{3}}.$$

From these facts, we can conclude that if

$$\sigma_{\min} \left( A_I^T A_I \right) > \|A\|_2^2 / 4$$

and the step size $\mu$ satisfies

$$\frac{1}{4\sigma_{\min} \left( A_I^T A_I \right)} < \mu < \|A\|_2^{-2},$$

then

$$0 < \lambda < 8e^{\frac{3}{2}} \sigma_{\min} \left( A_I^T A_I \right).$$

That is, under these conditions, $x^*$ must be a local minimizer of $L_{1/2}$ regularization. We state this formally as the following Theorem 2.

*Theorem 2:* Assume that the sequence $\{x^{(n)}\}$ generated by *half* algorithm converges to $x^*$. Let $I = supp(x^*)$. If

$$\sigma_{\min} \left( A_I^T A_I \right) > \frac{\|A\|_2^2}{4} \quad (15)$$

and

$$\frac{1}{4\sigma_{\min} \left( A_I^T A_I \right)} < \mu < \frac{1}{\|A\|_2^2}, \quad (16)$$

then $x^*$ is a local minimizer of $T_\lambda(x)$.

Note that the condition (15) implies $\frac{1}{4\sigma_{\min}(A_I^T A_I)} < \|A\|_2^{-2}$ naturally. (16) thus is a natural and reachable condition and, furthermore, whenever this condition is satisfied, the sequence $\{x^{(n)}\}$ is indeed convergent by Lemma 2. This shows that only the condition (15) is essential in Theorem 2. We notice that (15) is a concentration condition on singular values of the submatrix $A_I^T A_I$, and, in particular, it implies

$$\sigma_{\min} \left( A_I^T A_I \right) > \sigma_{\max} \left( A_I^T A_I \right) / 4,$$

or equivalently

$$Cond \left( A_I^T A_I \right) := \frac{\sigma_{\max} \left( A_I^T A_I \right)}{\sigma_{\min} \left( A_I^T A_I \right)} < 4, \quad (17)$$

where $Cond(A_I^T A_I)$ is the condition number of $A_I^T A_I$. (17) thus shows that the submatrix $A_I^T A_I$ is well-conditioned with the condition number lower than 4.

In recent years, a property called the restricted isometry property (RIP) of a matrix $A$ was introduced to characterize the concentration degree of the singular values of its submatrix with $k$ columns [3]. A matrix $A$ is said to be of the $k$-order RIP (denoted then by $\delta_k$-RIP) if there exists a $\delta_k \in (0,1)$ such that

$$(1 - \delta_k) \|x\|_2^2 \leq \|Ax\|_2^2 \leq (1 + \delta_k) \|x\|_2^2, \forall \|x\|_0 \leq k. \quad (18)$$

In other words, the RIP ensures that all submatrices of $A$ with $k$ columns are close to an isometry, and therefore distance-preserving. Let $k = \|x^*\|_0$. It can be seen from (18) that if $A$ possesses $\delta_k$-RIP with $\delta_k < 0.6$, then

$$Cond \left( A_I^T A_I \right) \leq \frac{1 + \delta_k}{1 - \delta_k} < 4.$$

Thus, we can claim that when $A$ satisfies a certain RIP, the condition (15) in Theorem 2 can be satisfied. In particular, we have the following proposition.

*Proposition 3:* Assume that $k < N/2$ and $A$ satisfies $\delta_k$-RIP with $\delta_k < \frac{3}{4+2N/k}$ or $\delta_{2k}$-RIP with $\delta_{2k} < \frac{3}{4+N/k}$, then (15) holds.

This can be directly checked by the facts that $\sigma_{\min}(A_I^T A_I) \geq 1 - \delta_k$, $\sigma_{\min}(A_I^T A_I) \geq 1 - \delta_{2k}$, $\sigma_{\max}(A^T A) \leq 1 + \delta_N$, $\delta_N \leq \frac{2N}{k}\delta_k$ and $\delta_N \leq \frac{N}{k}\delta_{2k}$ (c.f. Proposition 1 in [34]).

From Proposition 3, we can see, for instance, when $k/N = 1/3$ and $A$ satisfies $\delta_k$-RIP with $\delta_k < 3/10$ or $\delta_{2k}$-RIP with $\delta_{2k} < 3/7$, the condition (15) is satisfied, and therefore, by Theorem 2, the *half* algorithm converges to a local minimizer of $L_{1/2}$ regularization. It is noted that in the condition of Proposition 3, we always have $\delta_k < 3/8$ and $\delta_{2k} < 1/2$.

*Remark 3:* The RIP based conditions presented in Proposition 3 are just examples of sufficient conditions for the validity of Theorem 2. It is by no means the weakest one. Also, it is known that verifying the RIP for a given matrix $A$ is very often NP-hard [35]–[37], thus using RIP directly as a criterion of convergence of $L_{1/2}$ regularization is pessimistic in applications. On the positive side, however, in the field of compressed sensing, it has been demonstrated that if $A$ is generated via a certain random design, then $A$ satisfies the RIP conditions stated in Proposition 3 with high probability (c.f., [3], [6], [38], [39]). More specifically, if the entries of the matrix are chosen from a certain random distribution such as Gaussian, sub-Gaussian and Rademacher random distributions independently and identically, then the corresponding random matrix satisfies the restricted isometry property with high probability [3], [38]. Recently, [39] proved that if $A \in \mathbf{R}^{m \times N}$ is generated from a distribution satisfying the so-called concentration of measure inequality, then $A$ is of the RIP $\delta_k < \delta$ with the probability no less than $1 - 2e^{-C_2(\delta)m}$ provided that $m \geq C_1(\delta)k \log(N/k)$, where $C_1(\delta)$ and $C_2(\delta)$ are two constants only depend on $\delta$. Thus, in applications, if we have known that the measurement matrix $A$ was generated via certain a random distribution discussed above, we may verify the RIP conditions presented in Proposition 3 via comparing the number of measurements $m$ with the theoretical lower bound $\mathcal{O}(k \log(N/k))$, though the constant $C_1(\delta)$ may be very large in theory.

*Remark 4:* In the seminal paper [40], Candes and Tao considered the Dantzig selector method for recovering the sparse signal from the noisy measurements, which is concerned with solving the following convex optimization model

$$\min_{x \in \mathbf{R}^N} \|x\|_1 \ s.t. \ \left\|A^T(y - Ax)\right\|_\infty \leq \eta, \qquad (19)$$

where $\|\cdot\|_\infty$ denotes the $l_\infty$-norm, which provides the largest magnitude entry in a vector and $\eta$ is a positive constant parameter that controls the level of the noise. They showed that under a RIP related condition, the Dantzig selector can identify the true sparse signal stably with high probability. More specifically, the RIP related condition reads as

$$\delta_k + \theta_{k,2k} < 1, \qquad (20)$$

where $\theta_{k,2k}$ is the so called restricted orthogonality constant (ROC) of order $(k, 2k)$. For any positive integers $k_1$ and $k_2$, $\theta_{k_1,k_2}$ is defined to be the smallest non-negative number $\theta$ such that

$$|\langle Au, Av\rangle| \leq \theta \|u\|_2 \|v\|_2$$

for all $k_1$-sparse vector $u$ and $k_2$-sparse vector $v$ with disjoint supports. Recently, Cai and Zhang [41] improved the condition (21) to a sharper condition, i.e.,

$$\delta_k + \theta_{k,k} < 1. \qquad (21)$$

According to [41], the condition (21) is mostly weaker than $\delta_{2k} < 1/3$ and $\delta_{2k} < 1/2$, and therefore, weaker than the RIP conditions we formulated in Proposition 3. Nevertheless, it should be observed two remarkable differences between these two types of results: First, the conditions (20) and (21) are for the model not for a specific algorithm, while the condition (15) and Proposition 3 are, however, just for a specific algorithm (i.e., the *half* algorithm) not for $L_{1/2}$ regularization model itself. This difference makes usually not completely comparable between these two types of results. Second, the Dantzig selector model (19) is a convex model, while $L_{1/2}$ regularization model (6) is a nonconvex model.

*Remark 5:* For a nonconvex algorithm, it is generally hard to guarantee that the algorithm achieves a global minimum. However, there are still some feasible strategies that may be adopted to achieve this aim. In the following, we will give some discussions on this in the perspectives of both theory and practice. Theoretically, some stringent conditions such as mutual coherence [42] and restricted isometry property (RIP) [3] are often used in the field of compressed sensing to guarantee that a nonconvex algorithm can recover the sparsest signal, which is also a global minimizer of the optimization problem under such conditions. In practice, there may be several possible ways that can be adopted to a nonconvex algorithm such that the algorithm can achieve a global minimum indeed or in probability. One of the possible ways is to choose a good initialization. It is well known that the convergent point of a nonconvex algorithm is commonly dependent on the initialization. Once the initialization is taken so good that it lies in the neighborhood of the global minimizer, then the algorithm can possibly achieve such global minimizer. Therefore, in practice, we may take a solution of the corresponding $l_1$-minimization problem as an initialization for *half* algorithm, since the solution of the $l_1$-minimization can be usually regarded as a good approximation of the global minimizer. Another possible way is to adopt some strategies of the global optimization such as simulated annealing (SA) and genetic evolution (GE) to avoid achieving a local minimizer. Recently, Xu and Wang has adopted the simulated annealing (SA) procedure to the *half* algorithm and then proposed a new hybrid thresholding algorithm [43]. It was demonstrated that under certain assumptions, the hybrid thresholding algorithm can converge to a global minimum in probability.

## IV. EVENTUALLY LINEAR CONVERGENCE RATE

In Section III, we have justified the convergence of the *half* algorithm to a local minimizer of $L_{1/2}$ regularization. How fast the *half* algorithm converges then is a problem that needs to be further answered. In this section, we derive an upper bound for the rate of convergence. We mainly show that the *half* algorithm possesses an eventually linear convergence speed.

*Theorem 3:* Assume that the sequence $\{x^{(n)}\}$ generated by the *half* algorithm converges to $x^*$, $I = supp(x^*)$ and $e = \min_{i \in I}|x_i^*|$. If either of the following two conditions satisfies

(1) $\sigma_{\min}(A_I^T A_I) > 0$ and $0 < \lambda < 8e^{3/2}\sigma_{\min}(A_I^T A_I)$;

(2) $\sigma_{\min}(A_I^T A_I) > \frac{\|A\|_2^2}{4}$ and $\frac{1}{4\sigma_{\min}(A_I^T A_I)} < \mu < \|A\|_2^{-2}$,

then there exists a sufficiently large positive integer $n_0$ and a constant $\rho^* \in (0, 1)$ such that when $n > n_0$,

$$\left\|x^{(n+1)} - x^*\right\|_2 \leq \rho^* \left\|x^{(n)} - x^*\right\|_2.$$

We present the proof of Theorem 3 in Section VII.

By Theorems 1 and 2, under the conditions of Theorem 3, the convergent limit $x^*$ is also a local minimizer. Thus, in this setting, the *half* algorithm converges to a local minimizer of $L_{1/2}$ regularization with an eventually linear convergence speed.

According to the proof of Theorem 3, we can find that

$$\rho^* = \frac{2\rho}{1 + \sqrt{1 - 4\rho C^* \left(\frac{\lambda\mu}{2}\right)^{\frac{3}{2}}}}, \quad (22)$$

where

$$\rho = \frac{8e^{\frac{3}{2}}\left(1 - \mu\sigma_{\min}\left(A_I^T A_I\right)\right)}{8e^{\frac{3}{2}} - \lambda\mu}$$

and $C^* = \min\left\{(1-\rho)\left(\frac{\lambda\mu}{2}\right)^{-3/2}, \frac{1}{\sqrt{2}(\lambda\mu)^{3/2}\rho}\right\}$. Since $0 < \lambda < 8e^{3/2}\sigma_{\min}(A_I^T A_I)$, we have $0 < \rho < 1$, and, furthermore, through simplifying the formula (22), we can find

$$\rho^* < \frac{2\rho}{|2\rho - 1| + 1}. \quad (23)$$

which shows that the bound of $\rho^*$ can be monotonically increasing as $\rho$ increases. Thus, it can be directly checked that $\rho^* < 1$. Moreover, it can be seen from (22) also that $\rho$ is monotonically increasing with $\lambda$ and decreasing with $\sigma_{\min}(A_I^T A_I)$ and $\mu$. Consequently, when the regularization parameter $\lambda$ is taken smaller or $\sigma_{\min}(A_I^T A_I)$, $\mu$ taken larger, $\rho^*$ may be smaller and thus the *half* algorithm can converge faster. This latter observation coincides clearly with the common understanding that a lager step size will lead to a faster convergence. Accordingly, we recommend that in practical application of $L_{1/2}$ regularization, a larger step size $\mu$ will be taken.

## V. Comparisons With Related Work

In this section, we compare the obtained theoretical results on *half* algorithm with some other related algorithms and work.

The first class of the tightly related algorithms are the iterative shrinkage and thresholding (IST) methods for $L_q$ regularization. Among them, the *hard* algorithm and the *soft* algorithm are two representatives, which respectively solves the $L_0$ regularization and $L_1$ regularization [32], [33]. It was demonstrated in [32], [33] that when $\mu = 1$ both *hard* and *soft* algorithms can converge to a stationary point whenever $\|A\|_2 < 1$. These classical convergence results can be generalized when a step size parameter $\mu$ is incorporated with the IST procedures, and in this case, the convergence condition becomes

$$0 < \mu < \|A\|_2^{-2}. \quad (24)$$

It can be seen that (24) is the exact condition of the convergence of *half* algorithm we have imposed in Lemma 2, which then supports that the classical convergence results of IST has been extended to the nonconvex $L_{1/2}$ regularization case. Furthermore, it was shown in [31] that when the measurement matrix $A$ satisfies the so-called finite basis injective (FBI) property and the stationary point possesses a strict sparsity pattern, the *soft* algorithm can converge to a global minimizer of $L_1$ regularization with a linear convergence rate. Such result is not surprising because of the convexity of $L_1$ regularization. As for convergence speed of the *hard* algorithm, it was demonstrated in [32] that under the condition $\mu = 1$ and $\|A\|_2 < 1$, *hard* algorithm will converge to a local minimizer with an eventually linear convergence rate (as far as we know, no result was given however for the case when step size $\mu$ is taken into consideration). However, both as algorithms for solving nonconvex models, Theorem 3 reveals that the *half* algorithm shares the same eventual convergence speed with the *hard* algorithm.

A generic IST algorithm has been also developed for solving the following general nonconvex regularization model

$$\min_{x \in \mathbf{R}^N} \left\{F(x) + \sum_{i=1}^{N} \phi(|x_i|)\right\}, \quad (25)$$

where $F$ is assumed to be a proper lower-semicontinuous (l.s.c.) convex function with gradient $\nabla F(x)$ being Lipschitz continuous and $\phi$ a nonconvex regularization term satisfying certain assumptions. The generalized gradient projection (GGP) algorithm proposed in [45] and the general iterative shrinkage and thresholding (GIST) algorithm suggested in [46] are two typical algorithms of IST for solving the general model. Nevertheless, it was only justified that both GGP and GIST algorithms can converge subsequentially to a stationary point [45], [46] (that is, there is a subsequence of the algorithm that converges to a stationary point). However, as a specific case of GGP algorithm, we have justified that the *half* algorithm can assuredly converge to a local minimizer with an eventually linear convergence rate.

Besides IST, another class of tightly related algorithms are the reweighted techniques that have been also widely used for solution to $L_q$ regularization with $q \in (0, 1)$. Two well known examples of such reweighted techniques are the iteratively reweighted least squares (IRLS) method [47] and the reweighted $l_1$ minimization (IRL1) method [15]. When specified to $L_{1/2}$ regularization, the convergence analysis conducted in [48] shows that the IRLS method converges with an eventually superlinear convergence rate under the assumption $A$ possesses a certain nullspace property (NSP). Such superlinear convergence seems very attractive. It should be observed that this is actually at the cost of that at each step of IRLS, a least squares problem has to be solved, which is time-consuming. We will provide a simulation study in the next section to show that although the eventual superlinear convergence rate, IRLS performs much slower than the *half* algorithm in practice, even the *half* algorithm is known only with eventual linear convergence rate.

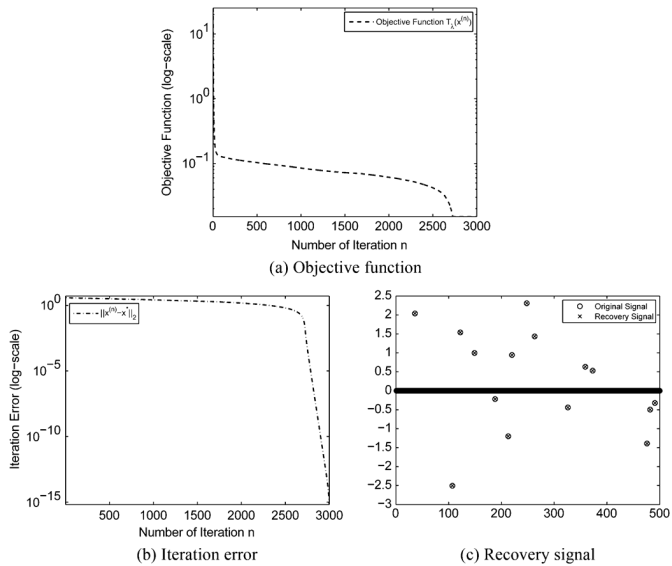

Fig. 1. Experiment for eventual linear convergence rate. (a) The trend of the objective function value, i.e., $T_\lambda(x^{(n)})$. (b) The trend of iteration error, i.e., $\|x^{(n)} - x^*\|_2$. (c) Recovery signal (Recovery MSE is $1.7928 \times 10^{-6}$). (a) Objective function; (b) iteration error; (c) recovery signal.

It was shown also in [49] that the IRL1 algorithm can converge to a stationary point and the eventual convergence speed is approximately linear. This convergence result seems to be worse than that we have obtained for the *half* algorithm (Theorem 3). Moreover, a series of simulations will be presented in the next section to support that the *half* algorithm can be faster than IRL1 algorithm in the perspective of time efficiency.

## VI. NUMERICAL EXPERIMENTS

We conduct a set of numerical experiments in this section to substantiate the validity of the theoretical analysis on convergence of the *half* algorithm.

### A. Convergence Rate Justification

We start with an experiment to confirm the linear rate of eventual convergence. For this purpose, given a sparse signal $x$ with dimension $N = 500$ and sparsity $k = 15$, shown as in Fig. 1(c), we considered the signal recovery problem through observation $y = Ax$, where the measurement matrix $A$ is of dimension $m \times N = 250 \times 500$ with Gaussian $\mathcal{N}(0, 1/250)$ i.i.d. entries. Such measurement matrix is known to satisfy (with high probability) the RIP with optimal bounds [38], [39] and thus the well-conditioned condition of Theorem 3 can be satisfied. We then applied the *half* algorithm to the problem with $\lambda = 0.001$ and $\mu = 0.99 \|A\|_2^{-2}$. The experiment results are reported in Fig. 1.

It can be seen from Fig. 1(a) that the iterative sequence of the objective function $\{T_\lambda(x^{(n)})\}$ is monotonically decreasing and converges to $T_\lambda^*$. Fig. 1(b) then shows how the recovery MSEs ($\|x^{(n)} - x^*\|_2$) of the iterations varies, from which we can see that, after approximately 2500 iterations, the *half* algorithm converges to a stationary point with a linear decay rate. Fig. 1(c) shows the comparison of recovered signal and the original signal, from which we see that the original sparse signal has been recovered with very high accuracy ($\|x^{(n)} - x^*\|_2 \leq 1.7928 \times 10^{-6}$). This experiment clearly justifies the convergence properties of the *half* algorithm we have verified in Lemma 2 and Theorems 3, particularly the expected eventually linear convergence rate of the *half* algorithm is substantiated.

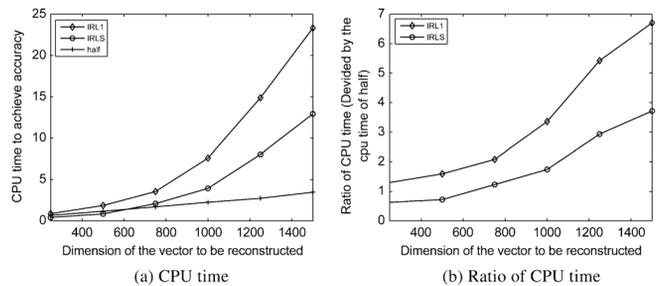

Fig. 2. Experiment for comparison of CPU times of different algorithms including *half*, IRLS and IRL1 algorithms. (a) The trends of CPU times of different algorithms. (b) The trends of the ratios of CPU times (divided by the cpu time of *half* algorithm). (a) CPU time; (b) ratio of CPU time.

### B. Comparisons With Reweighted Techniques

This set of experiments were conducted to compare the time costs of the *half* algorithm, IRLS algorithm [48] and IRL1 algorithm [15] for solving the same signal recovery problem with different settings $\{k, m, N\}$, where, as in Subsection 8.2 in [48], we took $k = 5$, $N = \{250, 500, 750, 1000, 1250, 1500\}$ and $m = N/5$. We implemented the three algorithms using Matlab without any specific optimization. In particular, we used the CVX Matlab package by Michael Grant and Stephen Boyd (http://www.stanford.edu/~boyd/cvx/) to perform the weighted $l_1$-minimization at each iteration step of IRL1 algorithm. Again, the measurement matrix $A$ was taken to be the $m \times N$ dimensional matrices with i.i.d. Gaussian $\mathcal{N}\left(0, \frac{1}{m}\right)$ entries. The experiment results are shown in Fig. 2.

As shown in Fig. 2(a), when $N$ is lower than 500, IRLS algorithm is slightly faster than the *half* algorithm. This is due to that in the low dimension cases, the computational burden of solving a low dimensional least square problem in IRLS is relatively low. Nevertheless, when $N > 500$, it can be observed that the *half* algorithm outperforms both IRLS and IRL1 algorithms in the perspective of CPU time. Furthermore, we can observe from Fig. 2(b) that as $N$ increases, the CPU times cost by IRL1 and IRLS algorithms increase much faster than the *half* algorithm, that is to say, the outperformance of the *half* algorithm in time cost can get more significant as dimension increases.

In addition, we analyze the computational complexity per iteration of three algorithms. At per iteration step of the *half* algorithm, some productions between matrix and vector are only required, and thus the computational complexity per iteration of the *half* algorithm is $\mathcal{O}(mN)$. While for IRLS algorithm, at each iteration, the solution of a least squares problem is required, of which the computational complexity is $\mathcal{O}(mN^2)$. Moreover, for IRL1 algorithm, at each iteration, a weighted $l_1$ minimization is required to solve, and thus the computational complexity per iteration of IRL1 algorithms is also $\mathcal{O}(mN^2)$. It demonstrates that when $N$ is large, then the computational complexities per iteration of IRLS and IRL1 algorithms are much higher than that of *half* algorithm.

## VII. PROOFS

In this section we present the proofs of Lemma 2, Lemma 3, Theorem 1 and Theorem 3.

### A. Proof of Lemma 2

The proof of Lemma 2 can follow that of Theorem 3 in [20]. We provide here, however, a much simplified proof. We first recall some known facts.

*Claim 1 ([20]):* Let $\{x^{(n)}\}$ be the sequence generated by the *half* algorithm with the step size $\mu$ satisfying $0 < \mu < \|A\|_2^{-2}$ and $T_\lambda(x) = \|Ax - y\|_2^2 + \lambda \|x\|_{1/2}^{1/2}$. Then
1) $\{x^{(n)}\}$ is a minimization sequence and $T_\lambda(x^{(n)})$ converges to $T_\lambda(x^*)$, where $x^*$ is any limit point of $\{x^{(n)}\}$;
2) $\{x^{(n)}\}$ is asymptotically regular, i.e., $\lim_{n\to\infty} \|x^{(n+1)} - x^{(n)}\|_2 = 0$.

By Claim 1, we easily conclude that the *half* algorithm is subsequentially convergent. We next claim that any limit point of $\{x^{(n)}\}$ is a stationary point of (7).

*Claim 2:* Under the setting of Claim 1, any limit point of $\{x^{(n)}\}$ is a stationary point of (7).

*Proof:* Denote $C_{\lambda,\mu}(x,z) = \|z - B_\mu(x)\|_2^2 + \lambda\mu\|z\|_{1/2}^{1/2}$ and let $D_{\lambda,\mu}(x) = C_{\lambda,\mu}(x,x) - \min_{z\in\mathbf{R}^N} C_{\lambda,\mu}(x,z)$. Then $D_{\lambda,\mu}(x) \geq 0$ and by Lemma 1, we have

$$D_{\lambda,\mu}(x) = 0 \text{ if and only if } x = F_{\lambda,\mu}(x).$$

Suppose that $x^*$ is a limit point $\{x^{(n)}\}$, then there is a subsequence $\{x^{(n)}\}$, say, $x^{(n_j)}$ such that $x^{(n_j)} \to x^*$ as $j \to \infty$. Since $x^{(n_j+1)} = F_{\lambda,\mu}(x^{(n_j)})$, we have

$$D_{\lambda,\mu}\left(x^{(n_j)}\right) = 2\left\langle \mu A^T\left(Ax^{(n_j)} - y\right), x^{(n_j)} - x^{(n_j+1)}\right\rangle$$
$$+ \lambda\mu\left(\left\|x^{(n_j)}\right\|_{\frac{1}{2}}^{\frac{1}{2}} - \left\|x^{(n_j+1)}\right\|_{\frac{1}{2}}^{\frac{1}{2}}\right) - \left\|x^{(n_j)} - x^{(n_j+1)}\right\|_2^2,$$

which implies

$$\lambda\left\|x^{(n_j)}\right\|_{\frac{1}{2}}^{\frac{1}{2}} - \lambda\left\|x^{(n_j+1)}\right\|_{\frac{1}{2}}^{\frac{1}{2}}$$
$$= \frac{1}{\mu}\left\|x^{(n_j)} - x^{(n_j+1)}\right\|_2^2 + \frac{1}{\mu}D_{\lambda,\mu}\left(x^{(n_j)}\right)$$
$$+ 2\left\langle A^T\left(Ax^{(n_j)} - y\right), x^{(n_j+1)} - x^{(n_j)}\right\rangle. \quad (26)$$

By (26), it then follows that

$$T_\lambda\left(x^{(n_j)}\right) - T_\lambda\left(x^{(n_j+1)}\right)$$
$$= \frac{1}{\mu}D_{\lambda,\mu}\left(x^{(n_j)}\right) + \frac{1}{\mu}\left\|x^{(n_j)} - x^{(n_j+1)}\right\|_2^2$$
$$- \left\|A\left(x^{(n_j)} - x^{(n_j+1)}\right)\right\|_2^2$$
$$\geq \frac{1}{\mu}D_{\lambda,\mu}\left(x^{(n_j)}\right) + (1/\mu - \|A\|_2^2)\left\|x^{(n_j)} - x^{(n_j+1)}\right\|_2^2.$$

Since $0 < \mu < \|A\|_2^{-2}$, we get

$$D_{\lambda,\mu}\left(x^{(n_j)}\right) \leq \mu\left(T_\lambda\left(x^{(n_j)}\right) - T_\lambda\left(x^{(n_j+1)}\right)\right).$$

and, by Claim 1, we have

$$T_\lambda\left(x^{(n_j)}\right) - T_\lambda\left(x^{(n_j+1)}\right) \to 0 \text{ as } j \to \infty,$$

This implies that $D_{\lambda,\mu}(x^*) = 0$, and hence $x^* = F_{\lambda,\mu}(x^*)$, as claimed. ∎

To complete the proof of Lemma 2, we need using also the following known property of the *half* algorithm.

*Claim 3 ([20]):* The number of the stationary points of (7) is finite.

With Claims 1–3, we can verify Lemma 2 as follows. From Claims 1 and 2, we first can conclude that the *half* algorithm converges subsequentially to a stationary point of (7). By Claim 3, any stationary point of (7) is isolated, which then, combined with the asymptotic regularity of the sequence (Claim 1(2)), implies that the whole sequence of $\{x^{(n)}\}$ must be convergent to a stationary point. Thus, the proof of Lemma 2 is completed.

### B. Proof of Lemma 3

By the assumption that $\{x^{(n)}\}$ converges to $x^*$, it is obvious that there is an $n_0$ such that whenever $n > n_0$, there holds $\|x^{(n)} - x^*\|_2 < \left(\frac{\lambda\mu}{2}\right)^{2/3}$.

We now prove by contradiction that the support set $I_n$ of $x^{(n)}$ is the same as $I$, the support set of $x^*$, whenever $n > n_0$. Assume this is not the case, namely, that $I_n \neq I$. Then we easily derive a contradiction through distinguishing the following two possible cases:

*Case 1:* $I_n \neq I$ and $I_n \cap I \subset I_n$. In this case, then there exists an $i_n$ such that $i_n \in I_n \setminus I$. By Proposition 1(2), it then implies

$$\left\|x^{(n)} - x^*\right\|_2 \geq \left|x_{i_n}^{(n)}\right| \geq \min_{i\in I_n}\left|x_i^{(n)}\right| > \left(\frac{\lambda\mu}{2}\right)^{\frac{2}{3}},$$

which contradicts to $\|x^{(n)} - x^*\|_2 < \left(\frac{\lambda\mu}{2}\right)^{2/3}$.

*Case 2:* $I_n \neq I$ and $I_n \cap I = I_n$. In this case, it is obvious that $I_n \subset I$. Thus, there exists an $i^*$ such that $i^* \in I \setminus I_n$. By Proposition 1(2), we still have

$$\left\|x^{(n)} - x^*\right\|_2 \geq |x_{i^*}^*| \geq \min_{i\in I}|x_i^*| > \left(\frac{\lambda\mu}{2}\right)^{\frac{2}{3}},$$

and it contradicts to $\|x^{(n)} - x^*\|_2 < \left(\frac{\lambda\mu}{2}\right)^{2/3}$.

Thus we have justified that $I_n = I$.

We further check $sign(x^{(n)}) = sign(x^*)$ by contradiction. As $I_n = I$ when $n > n_0$, it suffices to test that $sign(x_i^{(n)}) = sign(x_i^*)$ for any $i \in I$. Assume this is not the case. Then there exists an $i^* \in I$ such that $sign(x_i^{(n)}) \neq sign(x_i^*)$, and hence,

$$sign\left(x_{i^*}^{(n)}\right) sign\left(x_{i^*}^*\right) = -1.$$

From Proposition 1(2), it then implies

$$\left\|x^{(n)} - x^*\right\|_2 \geq \left|x_{i^*}^{(n)} - x_{i^*}^*\right| = \left|x_{i^*}^{(n)}\right| + |x_{i^*}^*|$$
$$\geq \min_{i\in I}\left\{\left|x_i^{(n)}\right| + |x_i^*|\right\} > 2\left(\frac{\lambda\mu}{2}\right)^{\frac{2}{3}},$$

contradicting again to $\|x^{(n)} - x^*\|_2 < \left(\frac{\lambda\mu}{2}\right)^{2/3}$. This contradiction shows $sign(x^{(n)}) = sign(x^*)$.

With this, the proof of Lemma 3 is completed.

## C. Proof of Theorem 1

Let $\beta^* = 1 - \left(\frac{\lambda}{8e^{3/2}\sigma_{\min}(A_I^T A_I)}\right)^{2/3}$. Since
$$0 < \lambda < 8e^{\frac{3}{2}}\sigma_{\min}\left(A_I^T A_I\right),$$
we have $0 < \beta^* < 1$. For any $\beta \in (0, 1)$, we define
$$c_\beta = \frac{1}{\left(\frac{\sqrt[3]{54}}{2(\lambda\mu)^{\frac{1}{3}}} + \frac{2}{\lambda}\|A_{I^c}^T A_I\|_2 \beta e\right)^2 e}.$$

Since $c_\beta$ is monotonically decreasing with respect to $\beta$, there exists a $\beta_0$ such that $\beta_0 \leq \beta^*$ and $\beta_0 \leq c_{\beta_0}$. We will proceed to show that whenever $\|h\|_2 < \beta_0 e$, it holds
$$T_\lambda(x^* + h) - T_\lambda(x^*) \geq 0.$$

Actually, we have
$$\begin{aligned}
T_\lambda(x^* &+ h) - T_\lambda(x^*) \\
&= 2h^T A^T(Ax^* - y) + h^T(A^T A)h \\
&\quad + \lambda \sum_{i=1}^N \left(|x_i^* + h_i|^{\frac{1}{2}} - |x_i^*|^{\frac{1}{2}}\right) \\
&= \lambda \sum_{i=1}^N \left(|x_i^* + h_i|^{\frac{1}{2}} - |x_i^*|^{\frac{1}{2}} + \frac{2}{\lambda}A_i^T(Ax^* - y)h_i\right) \\
&\quad + h^T(A^T A)h \\
&= \lambda \sum_{i \in I} \left(|x_i^* + h_i|^{\frac{1}{2}} - |x_i^*|^{\frac{1}{2}} + \frac{2}{\lambda}A_i^T(Ax^* - y)h_i\right) \\
&\quad + \lambda \sum_{i \in I^c}\left(|h_i|^{\frac{1}{2}} + \frac{2}{\lambda}A_i^T(Ax^* - y)h_i\right) + h^T(A^T A)h. \quad (27)
\end{aligned}$$

By Proposition 2, for any $i \in I^c$,
$$|h_i|^{\frac{1}{2}} + \frac{2}{\lambda}A_i^T(Ax^* - y)h_i \geq |h_i|^{\frac{1}{2}} - \frac{\sqrt[3]{54}}{2(\lambda\mu)^{\frac{1}{3}}}|h_i|, \quad (28)$$
and for any $i \in I$,
$$\begin{aligned}
|x_i^* + h_i|^{\frac{1}{2}} &- |x_i^*|^{\frac{1}{2}} + \frac{2}{\lambda}A_i^T(Ax^* - y)h_i \\
&= |x_i^* + h_i|^{\frac{1}{2}} - |x_i^*|^{\frac{1}{2}} - \frac{\operatorname{sign}(x_i^*)}{2|x_i^*|^{\frac{1}{2}}}h_i. \quad (29)
\end{aligned}$$

By Taylor expansion, for any $i \in I$, we then obtain
$$|x_i^* + h_i|^{\frac{1}{2}} - |x_i^*|^{\frac{1}{2}} - \frac{\operatorname{sign}(x_i^*)}{2|x_i^*|^{\frac{1}{2}}}h_i = \frac{-h_i^2}{8|x_i^* + \xi_i h_i|^{\frac{3}{2}}} \quad (30)$$

for some constant $\xi_i \in (0, 1)$. Furthermore, since
$$\max_{i \in I}|h_i| \leq \|h\|_2 < \beta_0 e,$$
we have for any $i \in I$,
$$|x_i^* + \xi_i h_i| \geq |x_i^*| - \xi_i|h_i| \geq (1 - \beta_0)e,$$
which shows
$$|x_i^* + h_i|^{\frac{1}{2}} - |x_i^*|^{\frac{1}{2}} + \frac{2}{\lambda}A_i^T(Ax^* - y)h_i \geq \frac{-h_i^2}{8(1-\beta_0)^{\frac{3}{2}}e^{\frac{3}{2}}}. \quad (31)$$

From (27)–(31), it thus follows that
$$\begin{aligned}
T_\lambda&(x^* + h) - T_\lambda(x^*) \\
&\geq \|Ah\|_2^2 + \lambda\sum_{i \in I^c}\left(|h_i|^{\frac{1}{2}} - \frac{\sqrt[3]{54}}{2(\lambda\mu)^{\frac{1}{3}}}|h_i|\right) \\
&\quad - \lambda\sum_{i \in I}\frac{h_i^2}{8(1-\beta_0)^{\frac{3}{2}}e^{\frac{3}{2}}} \\
&\geq \left(\|A_I h_I\|_2^2 - \frac{\lambda}{8(1-\beta_0)^{\frac{3}{2}}e^{\frac{3}{2}}}\|h_I\|_2^2\right) + \|A_{I^c}h_{I^c}\|_2^2 \\
&\quad + 2(A_I h_I)^T A_{I^c} h_{I^c} + \lambda\sum_{i \in I^c}\left(|h_i|^{\frac{1}{2}} - \frac{\sqrt[3]{54}}{2(\lambda\mu)^{\frac{1}{3}}}|h_i|\right) \\
&\geq \lambda\sum_{i \in I^c}\left(|h_i|^{\frac{1}{2}} - \frac{\sqrt[3]{54}}{2(\lambda\mu)^{\frac{1}{3}}}|h_i| - \frac{2}{\lambda}\left|A_i^T(A_I h_I)\right||h_i|\right) \\
&\quad + \left(\sigma_{\min}\left(A_I^T A_I\right) - \frac{\lambda}{8(1-\beta_0)^{\frac{3}{2}}e^{\frac{3}{2}}}\right)\|h_I\|_2^2. \quad (32)
\end{aligned}$$

Due to $\|h\|_2 < \beta_0 e \leq c_{\beta_0} e$, we find for any $i \in I^c$,
$$\begin{aligned}
|h_i|^{\frac{1}{2}} &\leq \frac{1}{\frac{\sqrt[3]{54}}{2(\lambda\mu)^{\frac{1}{3}}} + \frac{2}{\lambda}\|A_{I^c}^T A_I\|_2 \beta_0 e} \\
&\leq \frac{1}{\frac{\sqrt[3]{54}}{2(\lambda\mu)^{\frac{1}{3}}} + \frac{2}{\lambda}\|A_{I^c}^T A_I\|_2 \|h_I\|_2} \\
&\leq \frac{1}{\frac{\sqrt[3]{54}}{2(\lambda\mu)^{\frac{1}{3}}} + \frac{2}{\lambda}\|A_{I^c}^T A_I h_I\|_2} \\
&\leq \frac{1}{\frac{\sqrt[3]{54}}{2(\lambda\mu)^{\frac{1}{3}}} + \frac{2}{\lambda}\max_{i \in I^c}|A_i^T A_I h_I|} \\
&\leq \frac{1}{\frac{\sqrt[3]{54}}{2(\lambda\mu)^{\frac{1}{3}}} + \frac{2}{\lambda}|A_i^T A_I h_I|}. \quad (33)
\end{aligned}$$

Consequently,
$$\begin{aligned}
|h_i|^{\frac{1}{2}} &- \frac{\sqrt[3]{54}}{2(\lambda\mu)^{\frac{1}{3}}}|h_i| - \frac{2}{\lambda}\left|A_i^T A_I h_I\right||h_i| \\
&= \left(1 - \left(\frac{\sqrt[3]{54}}{2(\lambda\mu)^{\frac{1}{3}}} + \frac{2}{\lambda}\left|A_i^T(A_I h_I)\right|\right)|h_i|^{\frac{1}{2}}\right)|h_i|^{\frac{1}{2}} \\
&\geq 0. \quad (34)
\end{aligned}$$

By (32)–(34), we then obtain
$$T_\lambda(x^*+h) - T_\lambda(x^*) \geq \left(\sigma_{\min}\left(A_I^T A_I\right) - \frac{\lambda}{8(1-\beta_0)^{\frac{3}{2}}e^{\frac{3}{2}}}\right)\|h_I\|_2^2.$$

Since $\beta_0 \leq \beta^*$, this implies
$$\begin{aligned}
\sigma_{\min}&\left(A_I^T A_I\right) - \frac{\lambda}{8(1-\beta_0)^{\frac{3}{2}}e^{\frac{3}{2}}} \\
&\geq \sigma_{\min}\left(A_I^T A_I\right) - \frac{\lambda}{8(1-\beta^*)^{\frac{3}{2}}e^{\frac{3}{2}}} \\
&= 0.
\end{aligned}$$

and furthermore,

$$T_\lambda(x^* + h) - T_\lambda(x^*) \geq 0.$$

Thus, $x^*$ is a local minimizer of $L_{1/2}$ regularization. This implies Theorem 1.

### D. Proof of Theorem 3

Let

$$\rho = \frac{8e^{\frac{3}{2}}\left(1 - \mu\sigma_{\min}\left(A_I^T A_I\right)\right)}{8e^{\frac{3}{2}} - \lambda\mu}$$

and

$$C^* = \min\left\{(1-\rho)\left(\frac{2}{\lambda\mu}\right)^{\frac{3}{2}}, \frac{1}{\sqrt{2}(\lambda\mu)^{\frac{3}{2}}\rho}\right\}.$$

Since

$$0 < \lambda < 8e^{\frac{3}{2}}\sigma_{\min}\left(A_I^T A_I\right),$$

it follows that $0 < \rho < 1$ and $C^* > 0$. For any $0 < c < 1$, we define

$$\varepsilon_c = \max_{i \in I}\left\{\frac{3\lambda\mu|x_i^*|^{\frac{3}{2}}}{4\left(8|x_i^*|^{\frac{3}{2}} - \lambda\mu\right)\left(|x_i^*| - c\left(\frac{\lambda\mu}{2}\right)^{\frac{3}{2}}\right)^{\frac{5}{2}}}\right\}.$$

As $\varepsilon_c$ is monotonically increasing with respect to $c$, there exists a constant $c^*$ such that $0 < c^* < 1$ and $c^*\varepsilon_{c^*} < C^*$. By assumption, since $x^{(n)}$ converges to $x^*$, for a given constant $0 < c^* < 1$, there exists an $n_0$ such that when $n > n_0$,

$$\|x^{(n)} - x^*\|_2 < c^*\left(\frac{\lambda\mu}{2}\right)^{2/3}.$$

By Lemma 3, we have $sign(x^{(n)}) = sign(x^*)$ and $I_n = I$ when $n > n_0$. Denote

$$z^{(n)} = x^{(n)} - \mu A^T\left(Ax^{(n)} - y\right),$$
$$z^* = x^* - \mu A^T(Ax^* - y)$$

and $h^{(n)} = x^{(n)} - x^*$. Then $h_{I^c}^{(n)} = 0$ and $\|h^{(n)}\|_2 = \|h_I^{(n)}\|_2$ when $n > n_0$. For simplicity, we will assume henceforth $n > n_0$ is always satisfied.

On one hand, by Proposition 1(3), for any $i \in I$, $x_i^{(n+1)}$ and $x_i^*$ satisfy the following equations respectively,

$$x_i^{(n+1)} + \lambda\mu\frac{sign\left(x_i^{(n+1)}\right)}{4\sqrt{\left|x_i^{(n+1)}\right|}} = z_i^{(n)}, \quad (35)$$

$$x_i^* + \lambda\mu\frac{sign\left(x_i^*\right)}{4\sqrt{|x_i^*|}} = z_i^* \quad (36)$$

which, combined with $sign(x_i^{(n+1)}) = sign(x_i^*)$, then implies that for any $i \in I$,

$$x_i^{(n+1)} - x_i^* + \frac{\lambda\mu}{4}sign(x_i^*)\left(\frac{1}{\sqrt{\left|x_i^{(n+1)}\right|}} - \frac{1}{\sqrt{|x_i^*|}}\right) = z_i^{(n)} - z_i^*. \quad (37)$$

By applying Taylor's formula, we can conclude that for any $i \in I$, there exists a $\xi_{n+1} \in (0,1)$ such that

$$\frac{1}{\sqrt{\left|x_i^{(n+1)}\right|}} - \frac{1}{\sqrt{|x_i^*|}} = -\frac{sign(x_i^*)}{2|x_i^*|^{\frac{3}{2}}}h_i^{(n+1)} + \frac{3}{8\left|x_i^{\xi_{n+1}}\right|^{\frac{5}{2}}}\left|h_i^{(n+1)}\right|^2, \quad (38)$$

where $x_i^{\xi_{n+1}} = x_i^* + \xi_{n+1}(x_i^{(n+1)} - x_i^*)$. From (37) and (38), it then follows that

$$z_i^{(n)} - z_i^* = \left(1 - \frac{\lambda\mu}{8|x_i^*|^{\frac{3}{2}}}\right)h_i^{(n+1)} + \frac{3\lambda\mu sign(x_i^*)}{32\left|x_i^{\xi_{n+1}}\right|^{\frac{5}{2}}}\left|h_i^{(n+1)}\right|^2.$$

Denote by $\Lambda_1$ and $\Lambda_2$ the two diagonal matrices with

$$\Lambda_1(i,i) = 1 - \frac{\lambda\mu}{8|x_i^*|^{\frac{3}{2}}}$$

and

$$\Lambda_2(i,i) = \frac{3\lambda\mu sign(x_i^*)}{32\left|x_i^{\xi_{n+1}}\right|^{\frac{5}{2}}}.$$

Then we can write

$$z_I^{(n)} - z_I^* = \Lambda_1 h_I^{(n+1)} + \Lambda_2\left(h_I^{(n+1)} \odot h_I^{(n+1)}\right), \quad (39)$$

where $\odot$ is the Hadamard product or elementwise product.

On the other hand, for any $i \in I$, there holds the following identity

$$z_i^{(n)} - z_i^* = \left(x_i^{(n)} - x_i^*\right) - \mu\left(A_i^T\left(Ax^{(n)} - y\right) - A_i^T(Ax^* - y)\right).$$

Therefore, we have

$$z_I^{(n)} - z_I^* = \left(\mathbf{I} - \mu A_I^T A_I\right)h_I^{(n)}, \quad (40)$$

where $\mathbf{I}$ is the identity matrix.

Combining (39) and (40) gives

$$\Lambda_1 h_I^{(n+1)} + \Lambda_2\left(h_I^{(n+1)} \odot h_I^{(n+1)}\right) = \left(\mathbf{I} - \mu A_I^T A_I\right)h_I^{(n)}.$$

By applying Proposition 1(2), we easily check that for any $i \in I$, $|x_i^*| > \left(\frac{\lambda\mu}{2}\right)^{2/3}$, and

$$\Lambda_1(i,i) > 0.$$

Thus, it follows that

$$h_I^{(n+1)} = \Lambda_1^{-1} \left( \mathbf{I} - \mu A_I^T A_I \right) h_I^{(n)} - \Lambda_1^{-1} \Lambda_2 \left( h_I^{(n+1)} \odot h_I^{(n+1)} \right).$$

Furthermore, we have

$$\begin{aligned} \left\| h_I^{(n+1)} \right\|_2 &\leq \left\| \Lambda_1^{-1} \right\|_2 \left\| \mathbf{I} - \mu A_I^T A_I \right\|_2 \left\| h_I^{(n)} \right\|_2 \\ &\quad + \left\| \Lambda_1^{-1} \Lambda_2 \right\|_2 \left\| h_I^{(n+1)} \odot h_I^{(n+1)} \right\|_2 \\ &\leq \frac{8 e^{\frac{3}{2}} \left( 1 - \mu \sigma_{\min} \left( A_I^T A_I \right) \right)}{8 e^{\frac{3}{2}} - \lambda \mu} \left\| h_I^{(n)} \right\|_2 \\ &\quad + \left\| \Lambda_1^{-1} \Lambda_2 \right\|_2 \left\| h_I^{(n+1)} \odot h_I^{(n+1)} \right\|_2. \end{aligned}$$

Denote

$$\varepsilon_{n+1} = \max_{i \in I} \max_{x_i^{\xi_{n+1}}} \left\{ \left| \Lambda_1^{-1}(i,i) \Lambda_2(i,i) \right| \right\}.$$

Since

$$x_i^{\xi_{n+1}} = x_i^* + \xi_{n+1} \left( x_i^{(n+1)} - x_i^* \right)$$

and

$$\left| x_i^{(n+1)} - x_i^* \right| \leq \left\| x^{(n+1)} - x^* \right\|_2 \leq c^* \left( \frac{\lambda \mu}{2} \right)^{\frac{3}{2}},$$

we find

$$\varepsilon_{n+1} \leq \max_{i \in I} \left\{ \frac{3 \lambda \mu \left| x_i^* \right|^{\frac{3}{2}}}{4 \left( 8 \left| x_i^* \right|^{\frac{3}{2}} - \lambda \mu \right) \left( \left| x_i^* \right| - c^* \left( \frac{\lambda \mu}{2} \right)^{\frac{3}{2}} \right)^{\frac{5}{2}}} \right\}$$

$$= \varepsilon_{c^*}.$$

Consequently, we obtain

$$\left\| h_I^{(n+1)} \right\|_2 \leq \rho \left\| h_I^{(n)} \right\|_2 + \varepsilon_{c^*} \left\| h_I^{(n+1)} \right\|_2^2. \quad (41)$$

We further derive an relationship between $\|h_I^{(n)}\|_2$ and $\|h_I^{(n+1)}\|_2$ via the quadratic inequality (41). First, we observe that

$$1 - 4\rho \varepsilon_{c^*} \left\| h^{(n)} \right\|_2 > 1 - 4\rho \varepsilon_{c^*} c^* \left( \frac{\lambda \mu}{2} \right)^{\frac{2}{3}}$$

$$> 1 - 4\rho C^* \left( \frac{\lambda \mu}{2} \right)^{\frac{2}{3}} \geq (2\rho - 1)^2 \geq 0.$$

which implies

$$0 \leq \left\| h_I^{(n+1)} \right\|_2 \leq \frac{2\rho}{1 + \sqrt{1 - 4\rho \varepsilon_{c^*} \left\| h^{(n)} \right\|_2}} \left\| h^{(n)} \right\|_2.$$

Furthermore, since $\|h^{(n)}\|_2 < c^* \left( \frac{\lambda \mu}{2} \right)^{3/2}$, we conclude

$$\left\| h_I^{(n+1)} \right\|_2 < \frac{2\rho}{1 + \sqrt{1 - 4\rho \varepsilon_{c^*} c^* \left( \frac{\lambda \mu}{2} \right)^{\frac{3}{2}}}} \left\| h^{(n)} \right\|_2.$$

Denote

$$\rho^* = \frac{2\rho}{1 + \sqrt{1 - 4\rho \varepsilon_{c^*} c^* \left( \frac{\lambda \mu}{2} \right)^{\frac{3}{2}}}}.$$

Then, for sufficiently large $n$, we find

$$\left\| h_I^{(n+1)} \right\|_2 < \rho^* \left\| h_I^{(n)} \right\|_2.$$

Since

$$c^* \varepsilon_{c^*} < C^* \leq (1 - \rho) \left( \frac{\lambda \mu}{2} \right)^{-3/2},$$

it follows that

$$0 < \rho^* < 1.$$

Thus, the convergence of $\{x^{(n)}\}$ is eventually linear. This justifies Theorem 3 in the first condition case.

Note that the condition (2) implies the condition (1) in Theorem 3. Theorem 3 thus follows. This completes the proof of Theorem 3.

## VIII. CONCLUSION

We have conducted a study of a specific regularization framework, i.e., $L_{1/2}$ regularization, for better solution to the sparsity problem. Among various algorithms for $L_{1/2}$ regularization, *half* algorithm has been recognized one of the most effective algorithms. The main contribution of this paper is the establishment of the convergence analysis of *half* algorithm. In summary, we have verified that (i) the *half* algorithm converges to a stationary point in the loose circumstance $\mu \in (0, \|A\|_2^{-2})$; (ii) The algorithm can converge to a local minimizer of $L_{1/2}$ regularization if either the regularization parameter $\lambda$ is relatively small or the measurement matrix $A$ satisfies a certain concentration assumption (In particular, when $A$ possesses a certain RIP); (iii) The convergence speed is eventually linear. The obtained convergence result to stationary point generalizes those known for the *soft* and *hard* algorithms.

We have also provided a set of simulations to support the correctness of the established theoretical assertions. The efficiency of the *half* algorithm is further compared through simulations with the known reweighted techniques, another type of typical nonconvex regularization algorithms. The comparison shows that for higher dimensional problems and in terms of the time efficiency, the *half* algorithm outperforms the typical iteratively reweighted least squares (IRLS) algorithm and the iteratively reweighted $l_1$ minimization (IRL1) algorithm.

As shown in Fig. 1(c), given an appropriate regularization parameter, the *half* algorithm can indeed recover the true sparse signal with high precision. However, the corresponding theoretical guarantee have not been developed in the present research. We will investigate this problem in our future work.

**Jinshan Zeng** received the B.S. degree in Information and Computing Sciences from Xi'an Jiaotong University, Xi'an, China, in 2008. He is currently pursuing the Ph.D. degree with the School of Mathematics and Statistics, Xi'an Jiaotong University. He has been with the Department of Mathematics, University of California, Los Angeles, as a Visiting Scholar since Nov. 2013. His current research interests include sparse optimization, signal processing and synthetic aperture radar imaging.

**Shaobo Lin** received the BS degree in Mathematics and MS degree in Basic Mathematics from Hangzhou Normal University, China. Now he is a Ph.D. candidate in Xi'an Jiaotong University. His main research interests includes machine learning and scattered data fitting.

**Yao Wang** is currently working towards the Ph.D. degree in applied mathematics from Xi'an Jiaotong University, China. He was a research assistant at the Medical College of Georgia from Nov. 2009 to Oct. 2010, and a visiting scholar at Georgia Institute of Technology from Nov. 2012 to Nov. 2011. His research interests include computational biology, high-dimentsional statistical inference and statistical signal processing.

**Zongben Xu** received his Ph.D. degree in mathematics from Xi'an Jiaotong University, China, in 1987. He now serves as Vice President of Xi'an Jiaotong University, the Chief Scientist of National Basic Research Program of China (973 Project), and Director of the Institute for Information and System Sciences of the university. He is owner of the National Natural Science Award of China in 2007, and winner of CSIAM Su Buchin Applied Mathematics Prize in 2008. He delivered a 45 minute talk on the International Congress of Mathematicians 2010. He was elected as member of Chinese Academy of Science in 2011. His current research interests include intelligent information processing and applied mathematics.